\documentclass[12pt]{article}
\usepackage{amsthm,amsfonts,amsmath}
\usepackage{fullpage}

\newtheorem{theorem}{Theorem}

\newtheorem{lemma}[theorem]{Lemma}

\title{For each $\alpha > 2$ there is an infinite binary word with critical exponent $\alpha$}
\author{James D. Currie\thanks{The author's research was supported by an NSERC operating grant.} \&
Narad Rampersad\thanks{The author is supported by an NSERC Post-doctoral Fellowship.}\\
Department of Mathematics and Statistics\\
University of Winnipeg
\\
Winnipeg, Manitoba R3B 2E9\\CANADA\\
e-mail: {\tt currie@uwinnipeg.ca},\\
{\tt n.rampersad@uwinnipeg.ca}}

\begin{document}
\maketitle

\begin{abstract}
\noindent For each $\alpha > 2$ there is a binary word with critical exponent $\alpha$. \vspace{.1in}\\
\noindent Keywords: Combinatorics on words, repetitions,
critical exponent
\end{abstract}\maketitle

\small{ Mathematics Subject Classifications: 68R15}

\section{Introduction}
If $\alpha$ is a rational number, a word $w$ is an
$\alpha$~\emph{power} if there exist words $x$ and $x'$ and a
positive integer $n$, with $x'$ a prefix of $x$, such that $w =
x^nx'$ and $\alpha = n + |x'|/|x|$.  We refer to $|x|$ as a
\emph{period} of $w$. A word is $\alpha$~\emph{power-free} if none
of its subwords is a $\beta$~power with $\beta\ge \alpha$; otherwise,
we say the word \emph{contains an $\alpha$ power}.

The \emph{critical exponent} of an infinite word ${\bf w}$ is defined as
$$\text{sup} \lbrace \alpha \in \mathbb{Q} \mid {\bf w}
\text{ contains an $\alpha$-power} \rbrace.$$  Critical exponents of
certain classes of infinite words, such as Sturmian words
\cite{MP92,Van00} and words generated by iterated morphisms
\cite{Kri06,Kri07}, have received particular attention.

Krieger and Shallit \cite{KS07} proved that for every real number
$\alpha > 1$, there is an infinite word with critical exponent
$\alpha$.  As $\alpha$ tends to $1$, the number of letters required
to construct such words tends to infinity.  However, for $\alpha > 7/3$,
Shur \cite{Shu00} gave a construction over a binary alphabet.
For $\alpha > 2$, Krieger and Shallit gave a construction over a
four-letter alphabet and left it as an open problem to determine if
for every real number $\alpha \in (2,7/3]$, there is an infinite
binary word with critical exponent $\alpha$.  Currie, Rampersad,
and Shallit \cite{CRS06} gave examples of such words for a dense subset of real
numbers $\alpha$ in the interval $(2,7/3]$.  In this note we
resolve the question completely by demonstrating that for every real
number $\alpha > 2$, there is an infinite binary word with critical exponent
$\alpha$.

\section{Properties of the Thue-Morse morphism}
In this section we present some useful properties of the
\emph{Thue-Morse morphism}; \emph{i.e.}, the morphism $\mu$
defined by $\mu(0) = 01$ and $\mu(1) = 10$.

\begin{lemma}\label{2^s}
Let $s$ be a positive integer. Let $z$ be a subword of $\mu^s(01)$
with $|z|\ge 2^s$. Then $z$ does not have period $2^s$.
\end{lemma}
\begin{proof} Write $\mu^s(0)=a_1a_2\ldots a_{2^n}$, $\mu^s(1)=b_1b_2\ldots
b_{2^n}$. One checks by induction that $a_i=1-b_i$ for $1\le i\le
2^n$, and the result follows.
\end{proof}

%

Brandenburg \cite{Bra83} proved the following useful theorem, which
was independently rediscovered by Shur \cite{Shu00}.

\begin{theorem}[Brandenburg; Shur]
\label{shur}
Let $w$ be a binary word and let $\alpha > 2$ be a real number.  Then $w$
is $\alpha$~power-free if and only if $\mu(w)$ is $\alpha$~power-free.
\end{theorem}

The following sharper version of one direction of this theorem (implicit in
\cite{KS03}) is also useful.

\begin{theorem}\label{stronger}
Suppose $\mu(w)$ contains a subword $u$ of period $p$, with
$|u|/p>2$.  Then $w$ contains a subword $v$ of length
$\lceil|u|/2\rceil$ and period $p/2$.
\end{theorem}

We will also have call to use the deletion operator $\delta$ which
removes the first (left-most) letter of a word. For example,
$\delta(12345)=2345.$

\section{A binary word with critical exponent $\alpha$}
We denote by ${\cal L}$ the set of factors (subwords) of words of
$\mu(\{0,1\}^*).$

\begin{lemma}\label{r-1}
Let $00v\in{\cal L}$, and suppose that $00v$ is
$\alpha$~power-free, some fixed $\alpha>2$. Let
$r=\lceil\alpha\rceil$. Suppose that $0^{r}v=xuy$ where $u$ is an
$\alpha$~power. Then $x=\epsilon$ and $u=0^r$.
\end{lemma}
\begin{proof}
Suppose that $u$ has period $p$. Since $00v$ is
$\alpha$~power-free, we can write $u=0^sv'$, $x=0^{r-s}$ for some
integer $s$, $3\le s \le r$, some prefix $v'$ of $v$. If $0^p$ is
not a prefix of $u$ then the prefix of $u$ of length $p$ contains
the subword 0001. Since $\alpha>2$, this means that $0001$ is a
subword of $u$ at least twice, so that 0001 is a subword of $00v$.
This is impossible, since $00v\in{\cal L}$.

Therefore, $0^p$ is a prefix of $u$, and $u$ has the form $0^t$,
some integer $t\ge \alpha$. This implies that $u$ has $0^r$ as a
prefix, so that $x=\epsilon$ and $u=0^r$.
\end{proof}

\begin{lemma}\label{phi}
Let $\alpha>2$ be given, and let $r=\lceil\alpha\rceil$. Let $s,t$
be positive integers, such that $s\ge 3$ and there are words
$x,y\in\{0,1\}^*$ such that $\mu^s(0)=x00y$ with $|x|=t$. Suppose
that $2<r-t/2^s<\alpha$ and $00v\in{\cal L}$ is
$\alpha$~power-free. \begin{enumerate} \item{The word
$\delta^t\mu^s(0^rw)$ has a prefix which is a $\beta$~power, where
$\beta=r-t/2^s$.} \item{Suppose that $00v$ contains a
$\beta$~power of period $p$, some $\beta$ and $p$. Then
$\delta^t\mu^s(0^rw)$ contains a $\beta$~power of period
$2^sp$.}\item{Word $\delta^t\mu^s(0^rw)$ is
$\alpha$~power-free}\end{enumerate}
\end{lemma}
\begin{proof}
We start by observing that $\mu^s(0^r)$ has period $2^s$. It
follows that $\delta^t\mu^s(0^r)$ has period $2^s$, length
$r2^s-t$, and hence is a $(r2^s-t)/2^s=\beta$~power.

Now suppose $u$ is a $\beta$~power of period $p$ in $00v$. Then
$\mu(u)$ is a $\beta$~power of period $2^s$ in $\mu^s(00v)$.
However, $\mu^s(0^{r-1}v)$ is a suffix of $\delta^t\mu^s(0^rv)$,
since $t<2^s=|\mu^s(0)|$. Thus $\mu^s(u)$ is a $\beta$~power of
period $2^s$ in $\delta^t\mu^s(0^rv)$.

Next, note that $\mu^s(0^{r-1}v)$ does not contain any $\kappa$
power, $\kappa\ge \alpha$. Otherwise, by Theorem~\ref{stronger}
and induction, $0^{r-1}v$ contains a $\kappa$~power. This is
impossible by Lemma~\ref{r-1}.

Suppose then that $\delta^t\mu^s(0^rv)$ contains a $\kappa$~power
$\hat{u}$ of period $p$, $\kappa\ge\alpha$. Using induction and
Theorem~\ref{stronger}, $0^rv$ contains a $\kappa$~power $u$ of
period $p/2^s$. By Lemma~\ref{r-1}, the only possibility is
$u=0^r$, and $p/2^s=1$. Thus $p=2^s$.

Since $00v\in{\cal L}$, the first letter of $v$ is a 1. Since
$\hat{u}$ has period $2^s$, by Lemma~\ref{2^s} no subword of
$\mu^s(01)$ of length greater than $2^s$ occurs in $\hat{u}$. We
conclude that either $\hat{u}$ is a subword of
$\delta^t\mu^s(0^r)$, or of $\mu^s(v)$, and hence of
$\mu^s(0^{r-1}v)$. As this second case has been ruled out earlier,
we conclude that $|\hat{u}|\le |\delta^t\mu^s(0^r)|=r2^s-t.$ This
gives a contradiction: $\hat{u}$ is a $\kappa$ power, yet
$|\hat{u}|/p\le (r2^s-t)/2^s=\beta<\alpha.$
\end{proof}

By construction, $\delta^t\mu^s(0^rw)$ has the form $00\hat{v}$
where $00\hat{v}\in{\cal L}$.

We are now ready to prove our main theorem:

\begin{theorem} Let $\alpha>2$ be a real number. There is a word
over $\{0,1\}$ with critical exponent $\alpha$.
\end{theorem}
\begin{proof} Call a real number
$\beta<\alpha$ \emph{obtainable} if $\beta$ can be written
$\beta=r-t/2^s$, where $r,s,t$ are positive integers, $s\ge 3$,
and the word obtained by removing a prefix of length $t$ from
$\mu^s(0)$ begins with 00. We note that $\mu^3(0)=01101001$ and
$\mu^3(1)=10010110$ are of length 8, and both contain 00 as a
subword; for a given $s\ge 3$ it follows that $r$ and $t$ can be
chosen so that $\beta=r-t/2^s<\alpha$ and $|\alpha-\beta|\le
7/2^s$; by choosing large enough $s$, an obtainable number $\beta$
can be chosen arbitrarily close to $\alpha$.

Let $\{\beta_i\}$ be a sequence of obtainable numbers converging
to $\alpha$. For each $i$ write $\beta_i=r_i-t_i/2^{s_i}$, where
$r_i,s_i,t_i$ are positive integers, $s_i\ge 3$, and the word
obtained by removing a prefix of length $t_i$ from $\mu^{s_i}(0)$
begins with 00. If $00w\in{\cal L}$, denote by $\phi_i(w)$ the
word $\delta^{t_i}\mu^{s_i}(0^{r_i}w)$.

Consider the sequence of words
\begin{eqnarray*}
w_1&=&\phi_1(\epsilon)\\
w_2&=&\phi_1(\phi_2(\epsilon))\\
w_3&=&\phi_1(\phi_2(\phi_3(\epsilon)))\\
&\vdots&\\
w_n&=&\phi_1(\phi_2(\phi_3(\cdots(\phi_n(\epsilon))\cdots)))\\
&\vdots&\end{eqnarray*}

By the third part of Lemma~\ref{phi}, if $00w\in{\cal L}$ is
$\alpha$~power-free, then so is $\phi_i(w)$. Since $00\epsilon$ is
$\alpha$~power-free, each $w_i$ is therefore $\alpha$~power-free.

By the first and second parts of Lemma~\ref{phi}, $w_n$ contains
$\beta_i$~powers, $i=1,2,\ldots, n$.

Note that $\epsilon$ is a prefix of $\phi_{n+1}(\epsilon)$, so
that
$$w_n=\phi_1(\phi_2(\phi_3(\cdots(\phi_n(\epsilon))\cdots)))$$ is a
prefix of
$$\phi_1(\phi_2(\phi_3(\cdots(\phi_n(\phi_{n+1}(\epsilon)))\cdots)))=w_{n+1}.$$
We may therefore let $w=\lim_{n\rightarrow\infty}w_i$.

Since every prefix of $w$ is $\alpha$~power-free, $w$ is
$\alpha$~power-free but contains $\beta_i$ powers for each $i$.
The critical exponent of $w$ is therefore $\alpha$.
\end{proof}

The following question raised by Krieger and Shallit remains open:
for $\alpha > 1$, if $\alpha$ powers are avoidable on a $k$-letter alphabet,
does there exist an infinite word over $k$ letters with critical exponent
$\alpha$?  In particular, for $\alpha > \text{RT}(k)$, where $\text{RT}(k)$ denotes
the \emph{repetition threshold} on $k$ letters (see \cite{Car07}),
does there exist an infinite word over $k$ letters with critical exponent
$\alpha$?  We believe that the answer is ``yes''.

\end{document}